\documentclass[12pt,reqno]{amsart} 

\usepackage[all]{xy}

\def\theoremnumberstyle{plain} 
\def\proofname{Proof}
\usepackage{tom}

\allowdisplaybreaks[3]

\SelectTips{cm}{12}  
\CompileMatrices

\begin{document}

   \parindent0cm

   \title[Obstructed Families]{Some obstructed equisingular families of curves on
     surfaces in $\PC^3$}
   \author{Thomas Markwig}
   \address{Universit\"at Kaiserslautern\\
     Fachbereich Mathematik\\
     Erwin--Schr\"odinger--Stra\ss e\\
     D --- 67663 Kaiserslautern\\
     Tel. +496312052732\\
     Fax +496312054795
     }
   \email{keilen@mathematik.uni-kl.de}
   \urladdr{http://www.mathematik.uni-kl.de/\textasciitilde keilen}


   \date{September, 2003.}

     
   \begin{abstract}
     Very few examples of obstructed equsingular families of curves on
     surfaces other than $\PC^2$ are known.  Combining results from
     \cite{Wes04} and \cite{Hir92} with an idea from \cite{CC99} we give in the
     present paper series of examples of 
     families of irreducible curves with simple
     singularities on surfaces in $\PC^3$ which are
     not T--smooth, i.e. do not have the 
     expected dimension,  (Section
     \ref{sec:examples}) and we compare this with conditions (showing
     the same asymptotics)  which
     ensure the existence of a T--smooth component (Section
     \ref{sec:condition}). 
   \end{abstract}

   \maketitle

   Below we are going to construct two series of equisingular
   families of curves on surfaces in $\PC^3$. In both examples the families are
   obstructed in the sense that they do not have the expected
   dimension. However, while in the first example at least the existence of
   such curves was expected, the families in the second example were
   expected to be empty. It would be interesting to see if the
   equisingular families contain further components which are
   well--behaved. However, the families which we construct fail to
   satisfy the numerical conditions for the existence of such a
   component given in Section \ref{sec:condition} by a factor of
   two. We do not know whether the families are reducible or not,
   or if they are smooth.

   \section{Examples of obstructed families}\label{sec:examples}

     Throughout this section $\Sigma$ will denote a 
     smooth projective surface 
     in $\PC^3$ of degree $n\geq 2$, and $H$ will be a hyperplane
     section of $\Sigma$. $S=\{\ks_1,\ldots,\ks_s\}$ will be a finite
     set of \emph{simple}
     singularity types, that is the $\ks_i$ are of type $A_k$ (given by
     $x^2-y^{k+1}=0$, $k\geq 1$), $D_k$ (given by $x^2y-y^{k-1}=0$,
     $k\geq 4$), or $E_k$ (given by $x^3-y^4=0$, $x^3-xy^3=0$,
     or $x^3-y^5=0$ for $k=6,7,8$ respectively). 
     In general, for positive integers $r_1,\ldots,r_s$ and
     $d$  we denote by  $V_{|dH|}^{irr}(r_1\ks_1,\ldots,r_s\ks_s)$ 
     the family of irreducible curves in the linear system $|dH|$
     with precisely $r=r_1+\ldots+r_s$ singular points, $r_i$ of which
     are of the type $\ks_i$, $i=1,\ldots,s$, where $\ks_i$ may be any
     analytic type of an isolated singularity.
     $V_{|dH|}^{irr}(r_1\ks_1,\ldots,r_s\ks_s)$ is called
     \emph{T--smooth} or \emph{not obstructed} if it
     is smooth of the expected dimension
     \begin{multline*}
       \expdim\left(V_{|dH|}^{irr}(r_1\ks_1,\ldots,r_s\ks_s)\right)=\dim|dH|-\sum_{i=1}^sr_i\cdot
       \tau(\ks_i)\\=
       \frac{nd^2+(4n-n^2)d}{2}+\frac{n^3-6n^2+11n-6}{6}-\sum_{i=1}^s r_i\cdot\tau(\ks_i),
     \end{multline*}
     where $\tau(\ks)=\dim_\C \C\{x,y\}/\big\langle\frac{\partial
       f}{\partial x},\frac{\partial f}{\partial y},f\big\rangle$ is the
     Tjurina number of the singularity type $\ks$ given by the local
     equation $f=0$. Note that $\tau(A_k)=\tau(D_k)=\tau(E_k)=k$.

     In this note we give examples of such equisingular families of
     curves which are obstructed in the sense that they have dimension
     larger than the expected one. We use the idea by which
     Chiantini and Ciliberto in \cite{CC99} showed the existence of
     obstructed families of nodal curves. 

     Let us fix a plane $P$ in $\PC^3$,
     a point $p$ outside $P$, and a curve $C$ of degree $d>1$ in $P$. If we intersect
     the cone $K_{C,p}$ over $C$ with vertex $p$ with $\Sigma$, this gives a
     curve $C'=K_{C,p}\cap \Sigma$ in $|dH|$ which is determined by the choice of $C$ and
     $p$ (see Lemma \ref{lem:cone}). In particular, if $C$ varies in an $N$-dimensional family in
     $P$, then $C'$ varies in an $N$-dimensional family on
     $\Sigma$, and if $C$ is irreducible, then for a general choice of
     $p$ the curve $C'$ will be irreducible as well (see Lemma
     \ref{lem:kegel}). Moreover, if $C$ has a singular 
     point $q$ of (simple) singularity type $\ks$ and 
     $\Sigma$ meets the line joining $p$ and $q$ transversally in $n$
     points, then $C'$ will have a singularity of the same type in each
     of these points. 

     \begin{example}\label{ex:1}
       Fix the set $S=\{\ks_1,\ldots,\ks_s\}$ and let $m=\max\{\tau(\ks)\;|\;\ks\in S\}$. Suppose that
       $n>2m+4$ and $d>> n$, and let $r_1,\ldots,r_s\geq 0$ be such that
       \begin{displaymath}
         \frac{d^2+(4-n)d+2}{2}\;\;\;\leq \;\;\;\sum_{i=1}^s r_i\cdot\tau(\ks_i)
         \;\;\;\leq\;\;\; \frac{d^2+(4-n)d+2}{2}+m-1.
       \end{displaymath}
       Then 
       \begin{displaymath}
          \sum_{i=1}^s r_i\cdot\tau(\ks_i)\leq \frac{d^2+(4-n)d+2}{2}+m-1
         \leq \frac{d^2}{2}-m\cdot d-3.
       \end{displaymath}
       Hence, by \cite{Wes04} Remark 3.3.5 the family
       $V=V_d^{irr}(r_1\ks_1,\ldots,r_s\ks_s)$ of irreducible plane
       curves $C$ of degree $d$
       with precisely $r=r_1+\ldots+r_s$ singular points, $r_i$ of
       which are of type $\ks_i$, is
       non-empty, and we may estimate its dimension:
       \begin{multline*}
         \dim(V)\geq\expdim(V)=\frac{d(d+3)}{2}-\sum_{i=1}^sr_i\cdot\tau(\ks_i)\\
         \geq 
         \frac{d(d+3)}{2}-\frac{d^2+(4-n)d+2}{2}-m+1=\frac{n-1}{2}\cdot d-m.
       \end{multline*}
       By the above construction we see that hence the family of
       curves $C'$ satisfies
       \begin{displaymath}
         \dim\left(V_{|dH|}^{irr}(nr_1\ks_1,\ldots,nr_s\ks_s)\right)\geq \frac{n-1}{2}\cdot d-m.
       \end{displaymath}
       However, the expected dimension of this family is
       \begin{multline*}
         \expdim\left(V_{|dH|}^{irr}(nr_1\ks_1,\ldots,nr_s\ks_s)\right)\\
         =\frac{nd^2+(4n-n^2)d}{2}+\frac{n^3-6n^2+11n-6}{6}-\sum_{i=1}^s n\cdot r_i\cdot\tau(\ks_i)\\
         \leq \frac{nd^2+(4n-n^2)d}{2}+\frac{n^3-6n^2+11n-6}{6}-
          n\cdot\left(\frac{d^2+(4-n)d+2}{2}\right)\\
          = \frac{n^3-6n^2+5n-6}{6}.
       \end{multline*}
       For $d>>n$, more precisely for 
       \begin{displaymath}
         d>\frac{n^3-6n^2+5n-6+6m}{3n-3},
       \end{displaymath}
       the expected dimension will be smaller than the
       actual dimension, which proves that the family is obstructed.
       
       In particular, if $S=\{\ks\}$, $\ks\in\{A_k,D_k,E_k\}$, and 
       \begin{displaymath}
         r=\left\lceil \frac{d^2+(4-n)d+2}{2k}\right\rceil,
       \end{displaymath}
       then $V_{|dH|}(nr\ks)$ is obstructed, once $d>>n>3k+4$.
     \end{example}

     Note that in the previous example 
     \begin{displaymath}
       \expdim\left(V_{|dH|}(nr_1\ks_1,\ldots,nr_s\ks_s)\right)
       \geq
       \frac{n^3-6n^2+5n-6}{6}-n\cdot(m-1)>0,
     \end{displaymath}
     that is, the existence of curves in $|dH|$ with the given
     singularities was expected. This not so in the following
     example. 

     \begin{example}
       Let $k$ be an \emph{even}, positive integer, $m\geq 1$, $d=2(k+1)^m$, and
       \begin{displaymath}
         r=\frac{3\cdot(k+1)\cdot\big((k+1)^{2m}-1\big)}{(k+1)^2-1}.
       \end{displaymath}
       Hirano proved in \cite{Hir92} the existence of an irreducible
       plane curve of degree $d$ with precisely $r$ singular points all of type
       $A_k$. Thus the above construction shows that 
       \begin{displaymath}
         V_{|dH|}^{irr}(nrA_k)
       \end{displaymath}
       is non-empty. However, the expected dimension is
       \begin{multline*}
         \expdim\big(V_{|dH|}^{irr}(nrA_k)\big) =
         \frac{nd^2+(4n-n^2)d}{2}+\frac{n^3-6n^2+11n-6}{6}-
         knr\\
         =\left(2-\frac{3\cdot(k^2+k)}{k^2+2k}\right)\cdot (k+1)^{2m}
         +o\big((k+1)^m\big),
       \end{multline*}
       which is negative for $m$ sufficiently large, since
       \begin{displaymath}
         \frac{3\cdot(k^2+k)}{k^2+2k}>2.
       \end{displaymath}
       This shows that $V_{|dH|}^{irr}(nrA_k)$ is obstructed for
       sufficiently large $k$.
     \end{example}

     \section{Some remarks on conditions for T--smoothness}\label{sec:condition}

     Unless otherwise specified in this section $\Sigma$ will be an
     arbitrary smooth projective surface, $H$ a very ample divisor on
     $\Sigma$, and $\ks_1,\ldots,\ks_s$ arbitrary (not necessarily
     different) topological or analytical
     singularity types. As in Section \ref{sec:examples} we denote for $d\geq
     0$ by
     $V_{|dH|}^{irr}(\ks_1,\ldots,\ks_s)$ the equisingular family of
     irreducible curves in $|dH|$ with precisely $s$ singular points of
     types $\ks_1,\ldots,\ks_s$, and again the expected dimension is
     \begin{displaymath}
       \expdim\big(V_{|dH|}^{irr}(\ks_1,\ldots,\ks_s)\big)=
       \dim|dH|-\sum_{i=1}^s\tau(\ks_i).
     \end{displaymath}
     $V_{|dH|}^{irr}(\ks_1,\ldots,\ks_s)$ is called T--smooth if it is
     smooth of the expected dimension. 

     By \cite{Kei01} Theorem 1.2 and 2.3 (which is a slight
     improvement of \cite{KT02} Theorem 3.3 and Theorem 4.3) 
     there is a curve $C\in V_{|dH|}^{irr}(\ks_1,\ldots,\ks_s)$ if 
     \begin{itemize}
     \item $d\cdot H^2-g(H)\geq m_i+m_j$, and
     \item
       $h^1\big(\Sigma,\kj_{X(\underline{m};\underline{z})/\Sigma}((d-1)H)\big)=0$ 
       for $\underline{z}\in\Sigma^r$ very general,
     \end{itemize}
     where $\underline{m}=(m_1,\ldots,m_s)$ with $m_i=e^*(\ks_i)$, a
     certain invariant which only depends on $\ks_i$. Moreover,
     $V_{|dH|}^{irr}(\ks_1,\ldots,\ks_s)$ is T--smooth at this curve
     $C$ (see e.g. \cite{Shu99} Theorem 1). Finally, by \cite{AH00}
     Theorem 1.1 there is a number $d(m)$ depending only on
     $m=\max\{m_1,\ldots,m_s\}$, such that for all $d\geq d(m)$ and
     for $\underline{z}\in\Sigma^r$ very general the map
     \begin{displaymath}
       H^0\big(\Sigma,\ko_\Sigma((d-1)H)\big)\longrightarrow
       H^0\big(\Sigma,\ko_{X(\underline{m};\underline{z})/\Sigma}((d-1)H)\big) 
     \end{displaymath}
     has maximal rank. In particular, if 
     \begin{displaymath}
       \dim|(d-1)H|\geq\deg\big(X(\underline{m};\underline{z})\big)
       =\sum_{i=1}^s\frac{m_i\cdot(m_i+1)}{2},
     \end{displaymath}
     then
     $h^1\big(\Sigma,\kj_{X(\underline{m};\underline{z})/\Sigma}((d-1)H)\big)=0$. 
     This proves the following Proposition.

     \begin{proposition}
       Let $S=\{\ks_1,\ldots,\ks_s\}$ be a finite set of pairwise different
       topological or analytical singularity types. Then there exists
       a number $d(S)$ such that for all $d\geq d(S)$ and
       $r_1,\ldots,r_s\geq0$ satisfying
       \begin{equation}\label{eq:cond:1}
         \sum_{i=1}^s r_i\cdot
         \frac{e^*(\ks_i)\cdot\big(e^*(\ks_i)+1\big)}{2}<\dim|(d-1)H|
       \end{equation}
       the equisingular family $V_{|dH|}^{irr}(r_1\ks_1,\ldots,r_s\ks_s)$
       has a non-empty T--smooth component.  
     \end{proposition}

     In \cite{Shu01} upper bounds for $e^*(\ks)$ are given. 
     For a non-simple analytical singularity type we have
     \begin{displaymath}
       e^*(\ks)=e^a(\ks)\leq 3\sqrt{\mu(\ks)}-2
     \end{displaymath}
     where $\mu(\ks)$ is the Milnor number of $\ks$, and for any
     topological singularity type 
     \begin{displaymath}
       e^*(\ks)=e^s(\ks)\leq\frac{9}{\sqrt{6}}\cdot\sqrt{\delta(\ks)}-1,
     \end{displaymath}
     where $\delta(\ks)$ is the delta invariant of $\ks$.
 
     For \emph{simple} singularity types there are the better bounds
     \begin{center}
       \begin{tabular}{|l|c||l|c|}
         \hline
         $\ks $ & $ e^*(\ks)$& $\ks $ & $ e^*(\ks)$\\\hline\hline
         $A_1 $ & $                      2$&
         $D_4 $ & $                          3$\\
         $A_2 $ & $                      3$&
         $D_5 $ & $                          4$\\
         $A_k, k=3,\ldots,7 $ & $    4$&
         $D_k, k\leq 6,\ldots,10 $ & $   5$\\
         $A_k, k=8,\ldots,10 $ & $   5$&
         $D_k, k\leq 11,\ldots,13 $ & $  6$\\
         $A_k, k\geq 1 $ & $         \leq 2\cdot\left\lfloor\sqrt{k+5}\right\rfloor$&
         $D_k, k\geq 1 $ & $         \leq 2\cdot\left\lfloor\sqrt{k+7}\right\rfloor+1$\\
         $E_6 $ & $  4$&
         $E_7 $ & $  4$\\
         $E_8 $ & $  5$ &&
         \\\hline
       \end{tabular}
     \end{center}
     \medskip

     In particular, if $S=\{\ks_1,\ldots,\ks_s\}$ is a finite set of
     \emph{simple} singularities, then there is a $d(S)$ such that for
     all $d\geq d(S)$ and all $r_1,\ldots,r_s\geq 0$ satisfying
     \begin{equation}\label{eq:cond:2}
       2\cdot\sum_{i=1}^s r_i\cdot\Big(\tau(\ks_i)+o\big(\sqrt{\tau(\ks_i)}\big)\Big)\leq 
       \dim|dH|
     \end{equation}
     the family $V_{|dH|}^{irr}(nr_1\ks_1,\ldots,nr_s\ks_s)$ has a
     non-empty T--smooth component.

     The families in Example \ref{ex:1} fail to satisfy this
     condition roughly by the factor $2$. We thus cannot conclude that
     these families are reducible as we could in a similar situation
     in \cite{Kei04a}. 

     However, if we compare Condition \ref{eq:cond:1} respectively
     \ref{eq:cond:2} to the conditions in \cite{GLS00} or \cite{Kei04}
     which ensure that the equisingular family is T--smooth at
     \emph{every} point, the latter basically invole the square of the
     Tjurina number and are therefore much more restrictive. This, of
     course, was to be expected. 

     \section{Some remarks on cones}

     In this section we collect some basic properties on cones
     used for the construcion in Section \ref{sec:examples}, in
     particular the dimension counts.

     For points $p_1,\ldots,p_r\in \PC^3$ we will denote by
     $\overline{p_1\ldots p_r}$ the
     linear span in $\PC^3$ of  $p_1,\ldots,p_r$, i.e.\ the smallest
     linear subspace containing $p_1,\ldots,p_r$.

     Let $P\subset\PC^3$ be a plane, $C\subset P$ a curve, and $p\in
     \PC^3\setminus P$ a point. Then we denote by 
     \begin{displaymath}
       K_{C,p}=\bigcup_{q\in C}\overline{qp}
     \end{displaymath}
     the cone over $C$ with vertex $p$. Note that 
     \begin{displaymath}
       K_{C,p}=\bigcup_{q\in K_{C,p}}\overline{qp}
     \end{displaymath}
     and that
     \begin{displaymath}
       K_{C,p}\cap P=C.
     \end{displaymath}

     We first show that $C$ and $p$ fix the cone uniquely except when
     $C$ is a line.

     \begin{lemma}\label{lem:cone}
       Let $P\subset\PC^3$ be a plane, and
       $C\subseteq P$ be an irreducible curve which is not a line. 

       Then for $p,p'\in\PC^3$ with $p\not=p'$ we have that
       $K_{C,p}\not=K_{C,p'}$.
     \end{lemma}
     \begin{proof}
       Suppose there are points $p\not= p'$ such that
       $K_{C,p}=K_{C,p'}$. Choose a point $x\in
       C\setminus\overline{pp'}$ and let $E=\overline{xpp'}$. Then for
       any point $y\in\overline{xp}\subset K_{C,p}=K_{C,p'}$ we have
       \begin{displaymath}
         \overline{yp'}\subset K_{C,p'},
       \end{displaymath}
       and thus $E=\bigcup_{y\in\overline{xp}}\overline{yp'}\subset
       K_{C,p'}$. This, however, implies that the line
       \begin{displaymath}
         l=E\cap P \subseteq K_{C,p'}\cap P=C
       \end{displaymath}
       is contained in $C$, and since $C$ is irreducible we would have
       $C=l$ in contradiction to our assumption that $C$ is not a
       line. Hence, $K_{C,p}\not= K_{C,p'}$ for $p\not=p'$.
     \end{proof}

     Finally we show that for a general $p$ the cone $K_{C,p}$
     intersects $\Sigma$ in an irreducible curve.

     \begin{lemma}\label{lem:kegel}
       Let $\Sigma\subset\PC^3$ be a smooth projective surface,
       $P\subset\PC^3$ be a plane such that $P\not=\Sigma$, and
       $C\subseteq P$ an irreducible curve which is not a line and not
       contained in $\Sigma$. 
       Then for
       $p\in\PC^3\setminus P$ general $K_{C,p}\cap\Sigma$ is
       irreducible.
     \end{lemma}
     \begin{proof}
       Consider the linear system $\kl$ in $\PC^3$ which is given as
       the closure of 
       \begin{displaymath}
         \left\{K_{C,p}\;|\;p\in\PC^3\setminus P\right\},
       \end{displaymath}
       and set for $q\in\PC^3\setminus P$
       \begin{displaymath}
         \kl_q=\{D\in\kl\;|\;q\in D\}.
       \end{displaymath}
       
       First we show that for $q'\in C$ and $q\not\in P$ 
       \begin{equation}\label{eq:kegel:1}
         \bigcap_{p\in\overline{qq'}}K_{C,p}=C\cup\overline{qq'}.
       \end{equation}
       Choose pairwise different point
       $p_1,\ldots,p_n\in\overline{qq'}\setminus\{q,q'\}$. Suppose that
       there is a
       $z\in \bigcap_{i=1}^nK_{C,p_i}\setminus(C\cup
       \overline{qq'})$. Since $z\in K_{C,p_i}$ there is a unique
       intersection point
       \begin{displaymath}
         x_i=\overline{zp_i}\cap C,
       \end{displaymath}
       and these points $x_1,\ldots, x_n$ are pairwise different,
       since $z\not\in\overline{qq'}=\overline{p_ip_j}$ for
       $i\not=j$. However, 
       \begin{displaymath}
         x_i\in\overline{zp_i}\subset\overline{zp_ip_j}=\overline{zqq'}
       \end{displaymath}
       and $x_i\in C\subset P$, so that 
       \begin{displaymath}
         q',x_1,\ldots,x_n\in P\cap\overline{zqq'}
       \end{displaymath}
       and $q',x_1,\ldots,x_n$ are pairwise different collinear points
       on $C$. Since $C$ is irreducible but not a line, this implies
       $\deg(C)\geq n+1$. In particular, if $n\geq
       \deg(C)$, then
       \begin{displaymath}
         \bigcap_{i=1}^nK_{C,p_i}=C\cup\overline{qq'},
       \end{displaymath}
       which implies \eqref{eq:kegel:1}.

       Note that by \eqref{eq:kegel:1} for $q\in\PC^3\setminus P$
       \begin{displaymath}
         \bigcap_{D\in\kl_q}D\subseteq\bigcap_{K_{C,p}\in\kl_q}K_{C,p}
         =\bigcap_{q'\in C}\bigcap_{p\in\overline{qq'}}K_{C,p}
         =\bigcap_{q'\in C}\big(C\cup\overline{qq'}\big)
         =C\cup\{q\},
       \end{displaymath}
       and thus
       \begin{equation}\label{eq:kegel:2}
         \bigcap_{D\in\kl}D\subseteq\bigcap_{q\in\PC^3\setminus
           P}\bigcap_{D\in\kl_q}D=C.
       \end{equation}
       
       Consider now the linear systems
       \begin{displaymath}
         \kl_\Sigma=\{D\cap\Sigma\;|\;D\in\kl\}\;\;\;\mbox{ and }\;\;\; 
         \kl_{q,\Sigma}=\{D\cap\Sigma\;|\;D\in\kl_q\}=\{D\in\kl_\Sigma\;|\;q\in D\}.          
       \end{displaymath}
       
       Suppose that $\kl_\Sigma$ does not contain any irreducible
       curve. By \eqref{eq:kegel:2} and since $C\not\subset\Sigma$ the
       linear system $\kl_\Sigma$ has no fixed component. Thus by
       Bertini's Theorem $\kl_\Sigma$ must be composed with a pencil
       $\kb$, and since for a general point $q\in\Sigma$ the pencil $\kb$
       contains only one element, say $\widetilde{C}$, through $q$,
       the linear system $\kl_{q,\Sigma}$ has a fixed component
       $\widetilde{C}$. But then
       \begin{displaymath}
         \widetilde{C}\subseteq \bigcap_{D\in\kl_q}D\cap\Sigma=C\cap\Sigma.
       \end{displaymath}
       However, $C\cap\Sigma$ is zero-dimensional, while
       $\widetilde{C}$ has dimension one. 

       This proves that $\kl_\Sigma$ contains an irreducible element, and
       thus its general element is irreducible. In particular, for
       $p\in\PC^3\setminus P$ general $K_{C,p}\cap\Sigma$ is irreducible.      
     \end{proof}


\providecommand{\bysame}{\leavevmode\hbox to3em{\hrulefill}\thinspace}

\end{document}